\input amstex
\documentstyle{amsppt}
%----------------------------------------------------------------
% Title:     Perfect cuboids and multisymmetric polynomials.
% Author:    Ruslan Sharipov
% Comments:  AmSTeX, 12 pages, amsppt style
% MSC-class: 11D41, 11D72, 13A50, 13F20
%----------------------------------------------------------------
%           Replacement for output macro definition
%
\catcode`@=11
\redefine\output@{%
  \def\break{\penalty-\@M}\let\par\endgraf
  \ifodd\pageno\global\hoffset=105pt\else\global\hoffset=8pt\fi  
  \shipout\vbox{%
    \ifplain@
      \let\makeheadline\relax \let\makefootline\relax
    \else
      \iffirstpage@ \global\firstpage@false
        \let\rightheadline\frheadline
        \let\leftheadline\flheadline
      \else
        \ifrunheads@ %\let\makefootline\relax
        \else \let\makeheadline\relax
        \fi
      \fi
    \fi
    \makeheadline \pagebody \makefootline}%
  \advancepageno \ifnum\outputpenalty>-\@MM\else\dosupereject\fi
}
\def\Beta{\mathchar"0\hexnumber@\rmfam 42}
\catcode`\@=\active
%----------------------------------------------------------------
\nopagenumbers
\chardef\textvolna='176
\def\negskp{\hskip -2pt}

\def\Sym{\operatorname{Sym}}
\def\compos{\,\raise 1pt\hbox{$\sssize\circ$} \,}

\accentedsymbol\hatgamma{\kern 2pt\hat{\kern -2pt\gamma}}
\accentedsymbol\checkgamma{\kern 2.5pt\check{\kern -2.5pt\gamma}}
\def\blue#1{#1}

\catcode`#=11\def\diez{#}\catcode`#=6
\catcode`&=11\catcode`&=4
\catcode`_=11\def\podcherkivanie{_}\catcode`_=8
\catcode`~=11\def\volna{~}\catcode`~=\active
\def\mycite#1{\cite{\blue{#1}}\immediate\special{ps:
     ShrHPSdict begin /ShrBORDERthickness 0 def}}
\def\myciterange#1#2#3#4{\cite{\blue{#2#3#4}}\immediate\special{ps:
     ShrHPSdict begin /ShrBORDERthickness 0 def}}
\def\mytag#1{%
    \tag#1}
\def\mythetag#1{\thetag{\blue{#1}}\immediate\special{ps:
     ShrHPSdict begin /ShrBORDERthickness 0 def}}
\def\myrefno#1{\no#1}
\def\myhref#1#2{\blue{#2}\immediate\special{ps:
     ShrHPSdict begin /ShrBORDERthickness 0 def}}
\def\myEarXivlink{\myhref{http://arXiv.org}{http:/\negskp/arXiv.org}}

\def\mytheorem#1{\csname proclaim\endcsname{Theorem #1}}
\def\mytheoremwithtitle#1#2{\csname proclaim\endcsname{Theorem #1#2}}
\def\mythetheorem#1{\blue{#1}\immediate\special{ps:
     ShrHPSdict begin /ShrBORDERthickness 0 def}}
\def\mylemma#1{\csname proclaim\endcsname{Lemma #1}}
\def\mylemmawithtitle#1#2{\csname proclaim\endcsname{Lemma #1#2}}

\def\mycorollary#1{\csname proclaim\endcsname{Corollary #1}}

\def\mydefinition#1{\definition{Definition #1}}
\def\mythedefinition#1{\blue{#1}\immediate\special{ps:
     ShrHPSdict begin /ShrBORDERthickness 0 def}}
\def\myconjecture#1{\csname proclaim\endcsname{Conjecture #1}}
\def\myconjecturewithtitle#1#2{\csname proclaim\endcsname{Conjecture #1#2}}

%----------------------------------------------------------------
% Cyrillic fonts definition
\font\eightcyr=wncyr8

%----------------------------------------------------------------
\pagewidth{360pt}
\pageheight{606pt}
\topmatter
\title
Perfect cuboids and multisymmetric polynomials.
\endtitle
\author
Ruslan Sharipov
\endauthor
\address Bashkir State University, 32 Zaki Validi street, 450074 Ufa, Russia
\endaddress
\email\myhref{mailto:r-sharipov\@mail.ru}{r-sharipov\@mail.ru}
\endemail
\abstract
    A perfect Euler cuboid is a rectangular parallelepiped with integer edges
and integer face diagonals whose space diagonal is also integer. The problem 
of finding such parallelepipeds or proving their non-existence is an old unsolved 
mathematical problem. The Diophantine equations of a perfect Euler cuboid have
an explicit $S_3$ symmetry. In this paper the cuboid equations are factorized 
with respect to their $S_3$ symmetry in terms of multisymmetric polynomials.
Some factor equations are calculated explicitly. 
\endabstract
\endtopmatter
%\loadbold
%\loadeufb
\TagsOnRight
\document
% \input countstyle

%\special{header=resource.eps}
\head
1. Introduction.
\endhead
     The search for perfect cuboids has the long history since 1719 (see 
\myciterange{1}{1}{--}{39}). This history is presented as an adventure
story in \mycite{37}. Let $x_1$, $x_2$, $x_3$ be the edges of a cuboid and
let $d_1$, $d_2$, $d_3$ be its face diagonals. Then we have the equations
$$
\xalignat 2
&\hskip -2em
(x_1)^2+(x_2)^2-(d_3)^2=0,
&&(d_3)^2+(x_3)^2-L^2=0,\\
&\hskip -2em
(x_2)^2+(x_3)^2-(d_1)^2=0,
&&(d_1)^2+(x_1)^2-L^2=0,
\mytag{1.1}\\
&\hskip -2em
(x_3)^2+(x_1)^2-(d_2)^2=0,
&&(d_2)^2+(x_2)^2-L^2=0,
\endxalignat
$$
where $L$ is the space diagonal of the cuboid. In the case of a perfect
Euler cuboid the equations \mythetag{1.1} constitute a system of Diophantine
equations with respect to seven variables $x_1$, $x_2$, $x_3$, $d_1$, $d_2$, 
$d_3$, and $L$. In \mycite{40} the equations \mythetag{1.1} were reduced
to a single Diophantine equation with respect to four especially introduced
parameters $a$, $b$, $c$, and $u$. On the base of this equation in \mycite{41}
three cuboid conjectures were formulated. These conjectures are studied
(but not yet proved) in \myciterange{42}{42}{--}{44}.\par
     In the present paper we apply a quite different approach to the equations
\mythetag{1.1}. The equations \mythetag{1.1} possess a natural $S_3$ symmetry. 
Indeed, the symmetric group $S_3$ is composed by transformations of the set
of three numbers $\{1,2,3\}$:
$$
\hskip -2em
\sigma=\left(\matrix\format \c\kern 5pt &\kern 5pt\c\kern 5pt &\kern 5pt\c\\ 
\vspace{-0.1ex}
1 & 2 & 3\\ \downarrow & \downarrow & \downarrow\\
\vspace{-0.1ex}
\sigma 1 & \sigma 2 & \sigma 3 \endmatrix\right).
\mytag{1.2}
$$
The transformation \mythetag{1.2} is applied to the equations \mythetag{1.1}
as follows:
$$
\pagebreak 
\xalignat 3
&\hskip -2em
\sigma(x_i)=x_{\sigma i},
&&\sigma(d_i)=d_{\sigma i},
&&\sigma(L)=L.
\mytag{1.3}
\endxalignat
$$
Looking at \mythetag{1.1} and \mythetag{1.3}, one can easily see that the
system of equations \mythetag{1.1} in whole is invariant with respect to the 
transformations $\sigma\in S_3$.\par
     The main goal of this paper is to factorize the equations \mythetag{1.1}
with respect to the $S_3$ symmetry \mythetag{1.3}. We reach this goal by 
deriving some new equations from \mythetag{1.1}. These new equations are 
written in terms of the values of so-called multisymmetric polynomials (they
generalize well-known symmetric polynomials).
\par
\head
2. Multisymmetric polynomials.
\endhead
     Multisymmetric polynomials, which are also known as vector symmetric
polynomials, diagonally symmetric polynomials, McMahon polynomials etc, were
initially studied in \myciterange{45}{45}{--}{51} (see also later publications
\myciterange{52}{52}{--}{65}). Let's consider a set of variables arranged into
some $m\times n$ matrix as follows:
$$
\hskip -2em
M=\Vmatrix x_{11} &\hdots &x_{1n}\\
\vdots & \ddots & \vdots\\
x_{m1} &\hdots &x_{mn}
\endVmatrix
\mytag{2.1}
$$
The symmetric group $S_n$ acts upon the matrix \mythetag{2.1} by permuting its
columns:
$$
\hskip -2em
\sigma(x_{ij})=x_{i\kern 0.5pt\sigma\kern -1.2pt j}.
\mytag{2.2}
$$
\mydefinition{2.1} A polynomial $p\in\Bbb Q[x_{11},\ldots,x_{mn}]$ is called 
{\it multisymmetric\/} if it is invariant with respect to the action 
\mythetag{2.2} of the symmetric group $S_n$, i\.\,e\. if 
$$
p(x_{1\kern 0.5pt\sigma\kern -1.2pt 1},\,\ldots,x_{m\kern 0.5pt\sigma\kern -1.2pt
n})=p(x_{11},\,\ldots,x_{mn})\text{\ \ for all \ }\sigma\in S_n.
$$
\enddefinition
     Let $q(x)=q(x_{11},\,\ldots,x_{mn})$ be an arbitrary polynomial of the 
variables composing the matrix \mythetag{2.1}. Then we can produce a 
multisymmetric polynomial by applying the symmetrization operator $S$ 
to the polynomial $q(x_{11},\,\ldots,x_{mn})$:
$$
\hskip -2em
S(q(x_{11},\,\ldots,x_{mn}))=\frac{1}{n!}\sum_{\sigma\in S_n}\!q(x_{1\kern 0.5pt
\sigma\kern -1.2pt 1},\,\ldots,x_{m\kern 0.5pt\sigma\kern -1.2ptn}).
\mytag{2.3}
$$
Regular symmetric polynomials (see \mycite{66}) correspond to the special
case $m=1$ in the definition~\mythedefinition{2.1}. Like in the case $m=1$, in
general case $m>1$ there are elementary symmetric polynomials. However, in this
general case elementary symmetric polynomials are enumerated not by a single
index, but by a multiindex:
$$
\hskip -2em
\alpha=[\alpha_1,\,\ldots,\alpha_m]\text{, \ where \ }\alpha_i\geqslant 0
\text{\ \ and \ }|\alpha|=\alpha_1+\ldots+\alpha_m\leqslant n.
\mytag{2.4}
$$
Let's denote through $x^\alpha$ the following monomial:
$$
x^\alpha=\underbrace{x_{11}\cdot\ldots\cdot x_{1\kern 0.5pt\alpha_1}}_{\alpha_1}
\cdot\underbrace{x_{2\kern 0.5pt\alpha_1+1}\cdot\ldots\cdot x_{2\kern 0.5pt
\alpha_1+\alpha_2}}_{\alpha_2}\cdot\ldots\cdot\underbrace{x_{m\kern 0.5pt 
n-\alpha_m+1}\cdot\ldots\cdot x_{m\kern 0.5pt n}}_{\alpha_m}.\quad
\mytag{2.5}
$$
The variables in the product \mythetag{2.5} are taken from $n$ consecutive columns
of the matrix \mythetag{2.1}. The initial group of $\alpha_1$ of them is taken from 
the first row of this matrix, the next group of $\alpha_2$ of these variables is 
taken from the second row and so on. \pagebreak The last group of $\alpha_m$ variables 
is taken from the last $m$-th row of the matrix \mythetag{2.1}. If $\alpha_i=0$, then 
the corresponding $i$-th group in \mythetag{2.5} is empty and hence the variables of 
$i$-th row do not enter the monomial \mythetag{2.5} at all. 
\mydefinition{2.2} An elementary multisymmetric polynomial $e_\alpha(x_{11},\,\ldots,
x_{mn})$ corresponding to the multiindex \mythetag{2.4} is produced from the monomial 
\mythetag{2.5} by means of the symmetrization operator \mythetag{2.3} according to
the formula
$$
\hskip -2em
e_\alpha(x)=\frac{n!}{\alpha!}\,S(x^\alpha)\text{, \ where \ }
\alpha!=\alpha_1!\cdot\ldots\cdot\alpha_m!.
\mytag{2.6}
$$ 
\enddefinition
Note that the ratio $n!/\alpha!$ in \mythetag{2.6} is always an integer number
and $e_\alpha(x)$ is the sum of exactly $n!/\alpha!$ monomials produced from the 
monomial \mythetag{2.5} by means of the permutations of variables \mythetag{2.2}. 
In the case of the trivial multiindex $0=[0,\,\ldots,0]$ the formulas
\mythetag{2.5} and \mythetag{2.6} reduce to the following ones:
$$
\xalignat 2
&x^0=1,&&e_0=1.
\endxalignat
$$\par
    Like in the case of regular symmetric polynomials, there is the following 
fundamental theorem for multisymmetric polynomials. 
\mytheorem{2.1} The elementary multisymmetric polynomials \mythetag{2.6} with 
multiindi\-ces $0<|\alpha|\leqslant n$ generate the ring of all multisymmetric 
polynomials, i\.\,e\. each multisymmetric polynomial $p\in\Bbb Q[x_{11},\ldots,x_{mn}]$ 
can be expressed as a polynomial with rational coefficients through these elementary 
multisymmetric polynomials.
\endproclaim
    The proof of the fundamental theorem~\mythetheorem{2.1} can be found in 
\mycite{51}. Unfortunately the elementary multisymmetric polynomials \mythetag{2.6} 
are not algebraically independent over $\Bbb Q$ for $m>1$ (see \mycite{60}). For 
this reason the expression of $p$ as a polynomial with rational coefficients through 
the elementary multisymmetric polynomials $e_\alpha(x)$, which is claimed by the 
fundamental theorem~\mythetheorem{2.1}, is not unique.\par
\head
2. Multisymmetric polynomials associated with a cuboid.
\endhead
     Note that the formulas \mythetag{1.3} can be treated as a special case of the 
formulas \mythetag{2.2}. Indeed, let's compose the $2\times 3$ matrix
$$
\hskip -2em
M=\Vmatrix x_1 & x_2 &x_3\\
\vspace{1ex}
d_1 & d_2 & d_3\endVmatrix.
\mytag{3.1}
$$
Due to \mythetag{1.3} the transformations $\sigma\in S_3$ act as permutations 
of columns upon the matrix \mythetag{3.1}. Applying the 
definition~\mythedefinition{2.1} to the matrix \mythetag{3.1}, we get the 
concept of a multisymmetric polynomial of six variables $x_1,\,x_2,\,x_3$ and 
$d_1,\,d_2,\,d_3$. Now we calculate the elementary multisymmetric polynomials 
corresponding to the matrix \mythetag{3.1}. The first three of these polynomials 
are
$$
\align
&\hskip -2em
e_{\sssize [1,0]}=x_1+x_2+x_3,\\
&\hskip -2em
e_{\sssize [2,0]}=x_1\,x_2+x_2\,x_3+x_3\,x_1,
\mytag{3.2}\\
&\hskip -2em
e_{\sssize [3,0]}=x_1\,x_2\,x_3.
\endalign
$$
It is easy to see that the polynomials \mythetag{3.2} coincide with the regular 
symmetric polynomials of the three variables $x_1,\,x_2,\,x_3$. \pagebreak The next 
three elementary multisymmetric polynomials are similar to \mythetag{3.2}. They are
$$
\align
&\hskip -2em
e_{\sssize [0,1]}=d_1+d_2+d_3,\\
&\hskip -2em
e_{\sssize [0,2]}=d_1\,d_2+d_2\,d_3+d_3\,d_1,
\mytag{3.3}\\
&\hskip -2em
e_{\sssize [0,3]}=d_1\,d_2\,d_3.
\endalign
$$
The polynomials \mythetag{3.3} coincide with the regular symmetric polynomials 
of the three variables $d_1,\,d_2,\,d_3$. The rest of the elementary multisymmetric 
polynomials are actually multisymmetric. They include variables from both rows
of the matrix $M$: 
$$
\align
&\hskip -2em
e_{\sssize [2,1]}=x_1\,x_2\,d_3+x_2\,x_3\,d_1+x_3\,x_1\,d_2,\\
&\hskip -2em
e_{\sssize [1,1]}=x_1\,d_2+d_1\,x_2+x_2\,d_3+d_2\,x_3+x_3\,d_1+d_3\,x_1,
\mytag{3.4}\\
&\hskip -2em
e_{\sssize [1,2]}=x_1\,d_2\,d_3+x_2\,d_3\,d_1+x_3\,d_1\,d_2.
\endalign
$$
The polynomials \mythetag{3.2}, \mythetag{3.3}, and \mythetag{3.4} constitute the
complete set of elementary multisymmetric polynomials associated with the matrix
\mythetag{3.1}.\par
\head 
4. The first four factor equations.
\endhead
     Note that the variables $x_1,\,x_2,\,x_3$ and $d_1,\,d_2,\,d_3$ in the matrix
\mythetag{3.1} are not independent. They are related to each other by means of the 
polynomial equations \mythetag{1.1}. For this reason the elementary multisymmetric 
polynomials \mythetag{3.2}, \mythetag{3.3}, and \mythetag{3.4} produced from the 
variables $x_1,\,x_2,\,x_3$ and $d_1,\,d_2,\,d_3$ gain more algebraic relations 
in addition to those present in the case of independent variables $x_1,\,x_2,\,x_3$ 
and $d_1,\,d_2,\,d_3$ (see comments to fundamental theorem~\mythetheorem{2.1}).
These algebraic relations are written as polynomial equations with coefficients
in $\Bbb Q$:
$$
\hskip -2em
p(e_{\sssize [1,0]},e_{\sssize [2,0]},e_{\sssize [3,0]},e_{\sssize [0,1]},
e_{\sssize [0,2]},e_{\sssize [0,3]},e_{\sssize [2,1]},e_{\sssize [1,1]},
e_{\sssize [1,2]},L)=0.
\mytag{4.1}
$$
The polynomial equations of the form \mythetag{4.1} derived from \mythetag{1.1}
as well as those fulfilled identically due to \mythetag{3.2}, \mythetag{3.3},
and \mythetag{3.4} are called {\it factor equations\/} of the cuboid equations 
\mythetag{1.1} with respect to their $S_3$ symmetry. Our present goal is 
to reveal some of these factor equations explicitly.\par
     The equations \mythetag{1.1} are quadratic with respect to their variables.
For this reason it is quite likely that there are no {\bf linear relationships} 
between multisymmetric polynomials \mythetag{3.2}, \mythetag{3.3}, and \mythetag{3.4}. 
As for {\bf higher order relationships}, they do actually exist. In order to reveal
them we need to consider squares, cubes, fourth powers etc, and various mutual 
products of the multisymmetric polynomials \mythetag{3.2}, \mythetag{3.3}, and 
\mythetag{3.4}. For the square $(e_{\sssize [1,0]})^2$ we have 
$$
\hskip -2em
(e_{\sssize [1,0]})^2=x_1^2+x_2^2+x_3^2+2\,(x_1\,x_2+x_3\,x_1+x_3\,x_2).
\mytag{4.2}
$$
On the other hand, from the cuboid equations \mythetag{1.1} we derive
$$
\hskip -2em
x_1^2+x_2^2+x_3^2=L^2.
\mytag{4.3}
$$
Applying \mythetag{4.3} to \mythetag{4.2} and comparing the result with 
\mythetag{3.2}, we derive 
$$
\hskip -2em
(e_{\sssize [1,0]})^2-2\,e_{\sssize [2,0]}-L^2=0. 
\mytag{4.4}
$$
The equation \mythetag{4.4} is the first and the most simple factor equation
produced from the cuboid equations \mythetag{1.1}.\par
     The polynomial $e_{\sssize [0,1]}$ in \mythetag{3.3} is very similar to 
$e_{\sssize [1,0]}$. For its square we have
$$
\hskip -2em
(e_{\sssize [1,0]})^2=d_1^{\kern 0.9pt 2}+d_2^{\kern 0.9pt 2}+d_3^{\kern 0.9pt 2}
+2\,(d_1\,d_2+d_3\,d_1+d_3\,d_2).
\mytag{4.5}
$$
On the other hand, from the cuboid equations \mythetag{1.1} we derive
$$
\hskip -2em
d_1^{\kern 0.9pt 2}+d_2^{\kern 0.9pt 2}+d_3^{\kern 0.9pt 2}=2\,L^2.
\mytag{4.6}
$$
Applying \mythetag{4.6} to \mythetag{4.5} and comparing the result with 
\mythetag{3.3}, we derive 
$$
\hskip -2em
(e_{\sssize [0,1]})^2-2\,e_{\sssize [0,2]}-2\,L^2=0. 
\mytag{4.7}
$$
The equation \mythetag{4.7} is the second factor equation
produced from the cuboid equations \mythetag{1.1}. It is equally simple 
as the equation \mythetag{4.4}.\par
     In order to derive the third factor equation from the cuboid equations 
\mythetag{1.1} we consider the cube $(e_{\sssize [1,0]})^3$ and apply the
formula \mythetag{3.2}:
$$
\hskip -2em
\aligned
(e_{\sssize [1,0]})^3&=x_1^3+x_2^3+x_3^3+3\,x_1\,(x_2^2+x_3^2)\,+\\
&+\,3\,x_2\,(x_3^2+x_1^2)+3\,x_3(x_1^2+x_2^2)+6\,x_1\,x_2\,x_3.
\endaligned
\mytag{4.8}
$$
Using the cuboid equations \mythetag{1.1}, we derive the following formulas:
$$
\xalignat 2
&\hskip -2em
x_1^2=L^2-d_1^{\kern 0.9pt 2},
&&x_1^3=L^2\,x_1-d_1^{\kern 0.9pt 2}\,x_1,\\
&\hskip -2em
x_2^2=L^2-d_2^{\kern 0.9pt 2},
&&x_2^3=L^2\,x_2-d_2^{\kern 0.9pt 2}\,x_2,
\mytag{4.9}\\
&\hskip -2em
x_3^2=L^2-d_3^{\kern 0.9pt 2},
&&x_3^3=L^2\,x_3-d_3^{\kern 0.9pt 2}\,x_3.
\endxalignat
$$
Substituting \mythetag{4.9} into the equality \mythetag{4.8}, we obtain the formula
$$
\hskip -2em
\gathered
(e_{\sssize [1,0]})^3=-(x_1\,d_1^{\kern 0.9pt 2}+x_2\,d_2^{\kern 0.9pt 2}
+x_3\,d_3^{\kern 0.9pt 2})+7\,L^2\,(x_1+x_2+x_3)\,-\\
-\,3\,x_1\,(d_2^{\kern 0.9pt 2}+d_3^{\kern 0.9pt 2})
-3\,x_2\,(d_3^{\kern 0.9pt 2}+d_1^{\kern 0.9pt 2})
-3\,x_3\,(d_1^{\kern 0.9pt 2}+d_2^{\kern 0.9pt 2})+6\,x_1\,x_2\,x_3.
\endgathered
\mytag{4.10}
$$
The right hand side of the formula \mythetag{4.10} is a multisymmetric polynomial.
For this reason we can apply the theorem~\mythetheorem{2.1} to it. As a result we get
$$
\hskip -2em
\gathered
(e_{\sssize [1,0]})^3=2\,e_{\sssize [1,2]}+6\,e_{\sssize [3,0]}
+4\,e_{\sssize [0,2]}\,e_{\sssize [1,0]}\,-\\
-\,2\,e_{\sssize [0,1]}\,e_{\sssize [1,1]}-e_{\sssize [1,0]}
\,e_{\sssize [0,1]}^2+7\,e_{\sssize [1,0]}\,L^2.
\endgathered
\mytag{4.11}
$$
Note that the equation \mythetag{4.7} can be resolved with respect to
$e_{\sssize [0,2]}$:
$$
\hskip -2em
e_{\sssize [0,2]}=\frac{e_{\sssize [0,1]}^2}{2}-L^2.
\mytag{4.12}
$$
Applying \mythetag{4.12} to \mythetag{4.11}, we can write \mythetag{4.11} as follows:
$$
2\,e_{\sssize [1,2]}+6\,e_{\sssize [3,0]}-2\,e_{\sssize [0,1]}\,e_{\sssize [1,1]}
+e_{\sssize [1,0]}\,e_{\sssize [0,1]}^2+3\,e_{\sssize [1,0]}\,L^2-e_{\sssize [1,0]}^3=0.
\quad
\mytag{4.13}
$$
The equation \mythetag{4.13} is the third factor equation derived from the 
cuboid equations \mythetag{1.1}. It is more complicated than \mythetag{4.4} 
and \mythetag{4.7}.\par
     There is another way for transforming the cube $(e_{\sssize [1,0]})^3$ given by
the formula \mythetag{4.8}. Indeed, we can resolve the left column of the equations
\mythetag{1.1} with respect to $(x_1)^2$, $(x_2)^2$, and $(x_3)^2$. As a result we get
$$
\xalignat 3
&x_1^2=\frac{d_2^{\kern 0.9pt 2}+d_3^{\kern 0.9pt 2}-d_1^{\kern 0.9pt 2}}{2},
&&x_2^2=\frac{d_3^{\kern 0.9pt 2}+d_1^{\kern 0.9pt 2}-d_2^{\kern 0.9pt 2}}{2},
&&x_3^2=\frac{d_1^{\kern 0.9pt 2}+d_2^{\kern 0.9pt 2}-d_3^{\kern 0.9pt 2}}{2}.
\qquad\quad
\mytag{4.14}
\endxalignat
$$
The formulas \mythetag{4.14} can be used instead of the formulas in the left column 
of \mythetag{4.9}. Applying these formulas to \mythetag{4.8}, we can get an expression
analogous to \mythetag{4.10} and then we can continue transforming it in a way similar
to \mythetag{4.11} and \mythetag{4.12}, expecting to get some new equation similar to
\mythetag{4.13}. But actually we get the equation coinciding with \mythetag{4.13}.
\par
     Now let's consider the cube $(e_{\sssize [0,1]})^3$. It is given by the following
formula:
$$
\hskip -2em
\aligned
(e_{\sssize [0,1]})^3&=d_1^{\kern 0.9pt 3}+d_2^{\kern 0.9pt 3}+d_3^{\kern 0.9pt 3}
+3\,d_1\,(d_2^{\kern 0.9pt 2}+d_3^{\kern 0.9pt 2})\,+\\
&+\,3\,d_2\,(d_3^{\kern 0.9pt 2}+d_1^{\kern 0.9pt 2})
+3\,d_3\,(d_1^{\kern 0.9pt 2}+d_2^{\kern 0.9pt 2})+6\,d_1\,d_2\,d_3.
\endaligned
\mytag{4.15}
$$
Note that the equations of the left column of \mythetag{1.1} can be resolved with
respect to $(d_1)^2$, $(d_2)^2$, and $(d_3)^2$. They yield the equalities
$$
\xalignat 2
&\hskip -2em
d_1^{\kern 0.9pt 2}=x_2^2+x_3^2,
&&d_1^{\kern 0.9pt 3}=x_2^2\,d_1+x_3^2\,d_1,\\
&\hskip -2em
d_2^{\kern 0.9pt 2}=x_3^2+x_1^2,
&&d_2^{\kern 0.9pt 3}=x_3^2\,d_2+x_1^2\,d_2,
\mytag{4.16}\\
&\hskip -2em
d_3^{\kern 0.9pt 2}=x_1^2+x_2^2,
&&d_3^{\kern 0.9pt 3}=x_1^2\,d_3+x_2^2\,d_3.
\endxalignat
$$
Substituting \mythetag{4.16} into the equality \mythetag{4.15}, we obtain the formula
$$
\hskip -2em
\aligned
(e_{\sssize [0,1]})^3&=6\,d_1\,x_1^2+6\,d_2\,x_2^2+6\,d_3\,x_3^2+6\,d_1\,d_2\,d_3\,+\\
&+\,4\,d_1\,(x_2^2+x_3^2)+4\,d_2\,(x_3^2+x_1^2)+4\,d_3\,(x_1^2+x_2^2).
\endaligned
\mytag{4.17}
$$
The formula \mythetag{4.17} is analogous to the formula \mythetag{4.10}. Its right hand
side is a multisymmetric polynomial. For this reason we can apply the 
theorem~\mythetheorem{2.1} and get
$$
(e_{\sssize [0,1]})^3=2\,e_{\sssize [2, 1]}+6\,e_{\sssize [0,3]}
-2\,e_{\sssize [1,0]}\,e_{\sssize [1,1]}
-10\,e_{\sssize [2,0]}\,e_{\sssize [0,1]}
+6\,e_{\sssize [0,1]}\,e_{\sssize [1,0]}^2.
\mytag{4.18}
$$
Note that the equation \mythetag{4.4} can be resolved with respect to
$e_{\sssize [2,0]}$:
$$
\hskip -2em
e_{\sssize [2,0]}=\frac{1}{2}\,e_{\sssize [1,0]}^2-\frac{1}{2}\,L^2.
\mytag{4.19}
$$
Applying \mythetag{4.19} to \mythetag{4.18}, we can write the equality \mythetag{4.18} 
as follows:
$$
\hskip -2em
2\,e_{\sssize [2,1]}+6\,e_{\sssize [0,3]}-2\,e_{\sssize [1,0]}\,e_{\sssize [1,1]}
+e_{\sssize [0,1]}\,e_{\sssize [1,0]}^2+5\,e_{\sssize [0,1]}\,L^2-e_{\sssize [0,1]}^3=0.
\mytag{4.20}
$$
The equation \mythetag{4.20} is the fourth factor equation derived from the 
cuboid equations \mythetag{1.1}. \pagebreak It is similar to the equation 
\mythetag{4.13}.\par
     The equations of the second column in \mythetag{1.1} can also be resolved with 
respect to $(d_1)^2$, $(d_2)^2$, and $(d_3)^2$. Using them, we can write the formulas
$$
\xalignat 2
&\hskip -2em
d_1^{\kern 0.9pt 2}=L^2-x_1^2,
&&d_1^{\kern 0.9pt 3}=d_1\,L^2-d_1\,x_1^2,\\
&\hskip -2em
d_2^{\kern 0.9pt 2}=L^2-x_2^2,
&&d_2^{\kern 0.9pt 3}=d_2\,L^2-d_2\,x_2^2,
\mytag{4.21}\\
&\hskip -2em
d_3^{\kern 0.9pt 2}=L^2-x_3^2,
&&d_3^{\kern 0.9pt 3}=d_3\,L^2-d_3\,x_3^2.
\endxalignat
$$
The formulas \mythetag{4.21} can be used instead of the equations \mythetag{4.16}. As
a result we get another sequence of equations. However, the ultimate result appears
to be coinciding with the equation \mythetag{4.20}.\par
\head 
5. More factor equations.
\endhead
     In the next step we consider the square $(e_{\sssize [2,0]})^2$. Using the 
formulas \mythetag{3.2}, we get the following explicit expression for this square:
$$
(e_{\sssize [2,0]})^2=x_1^2\,x_2^2+x_2^2\,x_3^2+x_3^2\,x_1^2
+2\,x_1^2\,x_2\,x_3+2\,x_2^2\,x_3\,x_1+2\,x_3^2\,x_1\,x_2.\quad
\mytag{5.1}
$$
In order to transform \mythetag{5.1} we use the formulas \mythetag{4.9}. This yields
$$
\gathered
(e_{\sssize [2,0]})^2=d_1^{\kern 0.9pt 2}\,d_2^{\kern 0.9pt 2}
+d_2^{\kern 0.9pt 2}\,d_3^{\kern 0.9pt 2}+d_3^{\kern 0.9pt 2}\,d_1^{\kern 0.9pt 2}
-2\,L^2\,(d_1^2+d_2^2+d_3^2)-2\,(x_1\,x_2\,d_3^{\kern 0.9pt 2}\,+\\
+\,x_2\,x_3\,d_1^{\kern 0.9pt 2}+x_3\,x_1\,d_2^{\kern 0.9pt 2})
+2\,L^2\,(x_1\,x_2+x_2\,x_3+x_3\,x_1)+3\,L^4.
\endgathered\quad
\mytag{5.2}
$$
The right hand side of the formula \mythetag{5.2} is a multisymmetric polynomial.
For this reason we can apply the theorem~\mythetheorem{2.1} to it. As a result we get
$$
\hskip -2em
\gathered
(e_{\sssize [2,0]})^2=-2\,e_{\sssize [0,1]}\,e_{\sssize [0,3]}
+\frac{2}{3}\,e_{\sssize [1,0]}\,e_{\sssize [1,2]}-\frac{4}{3}\,e_{\sssize [0,1]}
\,e_{\sssize [2,1]}-\frac{2}{3}\,e_{\sssize [1,1]}^2\,+\\
+\,\frac{2}{3}\,e_{\sssize [0,1]}\,e_{\sssize [1,1]}\,e_{\sssize [1,0]}
+\frac{8}{3}\,e_{\sssize [2,0]}\,e_{\sssize [0,2]}
-\frac{2}{3}\,e_{\sssize [0,1]}^2\,e_{\sssize [2,0]}
+2\,e_{\sssize [2,0]}\,L^2\,-\\
-\frac{2}{3}\,e_{\sssize [1,0]}^2\,e_{\sssize [0,2]}
+e_{\sssize [0,2]}^2+4\,e_{\sssize [0,2]}\,L^2
-2\,e_{\sssize [0,1]}^2\,L^2+3\,L^4.
\endgathered
\mytag{5.3}
$$
Note that we can use the equation \mythetag{4.20} in order to express $e_{\sssize [0,3]}$
through the other elementary multisymmetric polynomials in \mythetag{4.20}:
$$
\hskip -2em
e_{\sssize [0,3]}=-\frac{1}{3}\,e_{\sssize [2,1]}
+\frac{1}{3}\,e_{\sssize [1,0]}\,e_{\sssize [1,1]}
+\frac{1}{6}\,e_{\sssize [0,1]}^3
-\frac{1}{6}\,e_{\sssize [0,1]}\,e_{\sssize [1,0]}^2
-\frac{5}{6}\,e_{\sssize [0,1]}\,L^2.
\mytag{5.4}
$$
Apart from \mythetag{5.4}, we apply the formulas \mythetag{4.12} and \mythetag{4.19} 
to \mythetag{5.3}. Then we get
$$
\hskip -2em
\aligned
8\,e_{\sssize [1,0]}\,&e_{\sssize [1,2]}-8\,e_{\sssize [0,1]}\,e_{\sssize [2,1]}
-8\,e_{\sssize [1,1]}^2+4\,e_{\sssize [0,1]}^2\,e_{\sssize [1,0]}^2\,-\\
&-\,e_{\sssize [0,1]}^4-3\,e_{\sssize [1,0]}^4
+10\,e_{\sssize [1,0]}^2\,L^2+4\,e_{\sssize [0,1]}^2\,L^2+L^4=0.
\endaligned
\mytag{5.5}
$$
The equation \mythetag{5.5} is the fifth factor equation derived from the 
cuboid equations \mythetag{1.1}. It is more complicated than all of the previous
factor equations.\par
     Now let's consider the other square $(e_{\sssize [0,2]})^2$. Using the formulas
\mythetag{3.3}, we get the following explicit expression for this square:
$$
(e_{\sssize [0,2]})^2=d_1^{\kern 0.9pt 2}\,d_2^{\kern 0.9pt 2}
+d_2^{\kern 0.9pt 2}\,d_3^{\kern 0.9pt 2}+d_3^{\kern 0.9pt 2}
\,d_1^{\kern 0.9pt 2}+2\,d_1^{\kern 0.9pt 2}\,d_2\,d_3
+2\,d_2^{\kern 0.9pt 2}\,d_3\,d_1+2\,d_3^{\kern 0.9pt 2}\,d_1\,d_2.
\quad
\mytag{5.6}
$$
In order to transform \mythetag{5.6} we use the formulas \mythetag{4.21}. This yields
$$
\gathered
(e_{\sssize [0,2]})^2=x_1^2\,x_2^2+x_2^2\,x_3^2+x_3^2\,x_1^2
-2\,L^2\,(x_1^2+x_2^2+x_3^2)-2\,(d_1\,d_2\,x_3^2\,+\\
+\,d_2\,d_3\,x_1^2+d_3\,d_1\,x_2^2)+2\,L^2\,(d_1\,d_2
+d_2\,d_3+d_3\,d_1)+3\,L^4.
\endgathered\quad
\mytag{5.7}
$$
The right hand side of the formula \mythetag{5.7} is a multisymmetric polynomial.
For this reason we can apply the theorem~\mythetheorem{2.1} to it. As a result we 
get
$$
\hskip -2em
\gathered
(e_{\sssize [0,2]})^2=
-2\,e_{\sssize [1,0]}\,e_{\sssize [3,0]}
-\frac{4}{3}\,e_{\sssize [1,0]}\,e_{\sssize [1,2]}
+\frac{2}{3}\,e_{\sssize [0,1]}\,e_{\sssize [2,1]}
-\frac{2}{3}\,e_{\sssize [1,1]}^2\,+\\
+\,\frac{2}{3}\,e_{\sssize [1,0]}\,e_{\sssize [0,1]}\,e_{\sssize [1,1]}
+\frac{8}{3}\,e_{\sssize [2,0]}\,e_{\sssize [0,2]}
-\frac{2}{3}\,e_{\sssize [0,1]}^2\,e_{\sssize [2,0]}
+4\,e_{\sssize [2,0]}\,L^2\,+\\
-\,\frac{2}{3}\,e_{\sssize [1,0]}^2\,e_{\sssize [0,2]}
+e_{\sssize [2,0]}^2+2\,e_{\sssize [0,2]}\,L^2
-2\,e_{\sssize [1,0]}^2\,L^2+3\,L^4.
\endgathered
\mytag{5.8}
$$
Note that we can use the equation \mythetag{4.13} in order to express 
$e_{\sssize [3,0]}$ through the other elementary multisymmetric polynomials 
in \mythetag{4.13}:
$$
\hskip -2em
e_{\sssize [3,0]}=-\frac{1}{3}\,e_{\sssize [1,2]}
+\frac{1}{3}\,e_{\sssize [0,1]}\,e_{\sssize [1,1]}
+\frac{1}{6}\,e_{\sssize [1,0]}^3
-\frac{1}{6}\,e_{\sssize [1,0]}\,e_{\sssize [0,1]}^2
-\frac{1}{2}\,e_{\sssize [1,0]}\,L^2.
\mytag{5.9}
$$
We apply the formulas \mythetag{5.9}, \mythetag{4.12}, and \mythetag{4.19} 
to \mythetag{5.8}. As a result we get
$$
\hskip -2em
\aligned
-8\,e_{\sssize [1,0]}\,&e_{\sssize [1,2]}+8\,e_{\sssize [0, 1]}\,e_{\sssize [2,1]}
-8\,e_{\sssize [1,1]}^2+4\,e_{\sssize [0,1]}^2\,e_{\sssize [1,0]}^2\,-\\
&-\,e_{\sssize [1,0]}^4-3\,e_{\sssize [0,1]}^4
+20\,e_{\sssize [0,1]}^2\,L^2-2\,e_{\sssize [1,0]}^2\,L^2-5\,L^4=0.
\endaligned
\mytag{5.10}
$$
The equation \mythetag{5.10} is the sixth factor equation derived from the 
cuboid equations \mythetag{1.1}. It is similar to the equation 
\mythetag{5.5}.\par
     The terms $e_{\sssize [1,2]}$, $e_{\sssize [2,1]}$, and $e_{\sssize [1,1]}$
are mixed in \mythetag{5.5} and \mythetag{5.10}. We can separate $e_{\sssize [1,2]}$
and $e_{\sssize [2,1]}$ from $e_{\sssize [1,1]}$ in the following two equations:
$$
\align
&\hskip -2em
\aligned
8\,e_{\sssize [1,0]}\,e_{\sssize [1,2]}&-8\,e_{\sssize [0,1]}\,e_{\sssize [2,1]}
+e_{\sssize [0,1]}^4-e_{\sssize [1, 0]}^4\,-\\
&-\,8\,e_{\sssize [0,1]}^2\,L^2+6\,e_{\sssize [1,0]}^2\,L^2+3\,L^4=0.
\endaligned
\mytag{5.11}\\
\vspace{1ex}
&\hskip -2em
\aligned
4\,e_{\sssize [1,1]}^2-2\,e_{\sssize [0,1]}^2\,&e_{\sssize [1,0]}^2
+e_{\sssize [0,1]}^4+e_{\sssize [1,0]}^4\,-\\
&\kern -1em -\,6\,e_{\sssize [0,1]}^2\,L^2-2\,e_{\sssize [1,0]}^2\,L^2+L^4=0.
\endaligned
\mytag{5.12}
\endalign
$$
The above equations \mythetag{5.11} and \mythetag{5.12} are just linear 
combinations of the fifth and the sixth factor equations \mythetag{5.5} 
and \mythetag{5.10}.\par
     In order to derive the next factor equation we consider the product
$e_{\sssize [2,0]}\,e_{\sssize [3,0]}$. Using the formulas
\mythetag{3.2}, we get the following explicit expression for this product:
$$
\hskip -2em
e_{\sssize [2,0]}\,e_{\sssize [3,0]}
=x_1\,x_2^2\,x_3^2+x_2\,x_3^2\,x_1^2+x_3\,x_1^2\,x_2^2.
\mytag{5.13}
$$
In order to transform \mythetag{5.13} we use the formulas \mythetag{4.9}. 
This yields
$$
\hskip -2em
\gathered
e_{\sssize [2,0]}\,e_{\sssize [3,0]}=x_1\,d_2^{\kern 0.9pt 2}
\,d_3^{\kern 0.9pt 2}+x_2\,d_3^{\kern 0.9pt 2}\,d_1^{\kern 0.9pt 2}
+x_3\,d_1^{\kern 0.9pt 2}\,d_2^{\kern 0.9pt 2}
-L^2\,x_1\,(d_2^{\kern 0.9pt 2}+d_3^{\kern 0.9pt 2})\,-\\
-\,L^2\,x_2\,(d_3^{\kern 0.9pt 2}+d_1^{\kern 0.9pt 2})
-L^2\,x_3\,(d_1^{\kern 0.9pt 2}+d_2^{\kern 0.9pt 2})
+L^4\,(x_1+x_2+x_3).
\endgathered
\mytag{5.14}
$$
The right hand side of the formula \mythetag{5.14} is a multisymmetric polynomial.
For this reason it can be expressed through elementary multisymmetric polynomials:
$$
\hskip -2em
\gathered
e_{\sssize [2,0]}\,e_{\sssize [3,0]}=
-e_{\sssize [1,1]}\,e_{\sssize [0,3]}+e_{\sssize [1,2]}\,e_{\sssize [0,2]}
+e_{\sssize [1,2]}\,L^2\,+\\
+\,e_{\sssize [0,2]}\,e_{\sssize [1,0]}\,L^2
-e_{\sssize [1,1]}\,e_{\sssize [0,1]}\,L^2
+e_{\sssize [1,0]}\,L^4.
\endgathered
\mytag{5.15}
$$
Now we transform \mythetag{5.15} with the use of the formulas \mythetag{5.9}, 
\mythetag{5.4}, \mythetag{4.19}, and \mythetag{4.12}: 
$$
\hskip -2em
\gathered
-4\,e_{\sssize [1,0]}\,e_{\sssize [1,1]}^2
+4\,e_{\sssize [1,1]}\,e_{\sssize [2,1]}
-2\,e_{\sssize [1,1]}\,e_{\sssize [0,1]}^3
+6\,e_{\sssize [1,2]}\,e_{\sssize [0,1]}^2\,+
\\
+\,2\,e_{\sssize [1,2]}\,e_{\sssize [1,0]}^2
+e_{\sssize [1,0]}^3\,e_{\sssize [0,1]}^2
-e_{\sssize [1,0]}^5
-2\,e_{\sssize [1,2]}\,L^2\,+
\\
+\,5\,e_{\sssize [1,0]}\,e_{\sssize [0,1]}^2\,L^2
+4\,e_{\sssize [1,0]}^3\,L^2-3\,e_{\sssize [1,0]}\,L^4=0.
\endgathered
\mytag{5.16}
$$
Then we apply \mythetag{5.12} to \mythetag{5.16} in order to eliminate the
term with the square $e_{\sssize [1,1]}^2$:
$$
\hskip -2em
\gathered
4\,e_{\sssize [1,1]}\,e_{\sssize [2,1]}
-2\,e_{\sssize [1,1]}\,e_{\sssize [0,1]}^3
+6\,e_{\sssize [1,2]}\,e_{\sssize [0,1]}^2\,+\\
+\,2\,e_{\sssize [1,2]}\,e_{\sssize [1,0]}^2
-\,e_{\sssize [1,0]}^3\,e_{\sssize [0,1]}^2
+e_{\sssize [1,0]}\,e_{\sssize [0,1]}^4
-2\,e_{\sssize [1,2]}\,L^2\,-\\
-\,e_{\sssize [1,0]}\,e_{\sssize [0,1]}^2\,L^2
+2\,e_{\sssize [1,0]}^3\,L^2
-2\,e_{\sssize [1,0]}\,L^4=0.
\endgathered
\mytag{5.17}
$$
The equation \mythetag{5.17} is the seventh factor equation derived from the 
cuboid equations \mythetag{1.1}. Its order is higher than the order of the
equations \mythetag{5.5} and \mythetag{5.10}.\par
     The eighth factor equation is derived similarly. In order to derive it
we consider the product $e_{\sssize [0,2]}\,e_{\sssize [0,3]}$. Applying 
the formulas \mythetag{3.3} to this product, we get
$$
\hskip -2em
e_{\sssize [0,2]}\,e_{\sssize [0,3]}=
d_1\,d_2^{\kern 0.9pt 2}\,d_3^{\kern 0.9pt 2}
+d_2\,d_3^{\kern 0.9pt 2}\,d_1^{\kern 0.9pt 2} 
+d_3\,d_1^{\kern 0.9pt 2}\,d_2^{\kern 0.9pt 2}.
\mytag{5.18}
$$
In order to transform \mythetag{5.18} we use the formulas \mythetag{4.21}. 
This yields
$$
\hskip -2em
\gathered
e_{\sssize [0,2]}\,e_{\sssize [0,3]}=d_1\,x_2^2\,x_3^2
+d_2\,x_3^2\,x_1^2+d_3\,x_1^2\,x_2^2
-L^2\,d_1\,(x_2^2+x_3^2)\,-\\
-\,L^2\,d_2\,(x_3^2+x_1^2)-L^2\,d_3\,(x_1^2+x_2^2)
+L^4\,(d_1+d_2+d_3).
\endgathered
\mytag{5.19}
$$
The right hand side of the formula \mythetag{5.19} is a multisymmetric polynomial.
For this reason it can be expressed through elementary multisymmetric polynomials:
$$
\hskip -2em
\gathered
e_{\sssize [0,2]}\,e_{\sssize [0,3]}=
-e_{\sssize [1,1]}\,e_{\sssize [3,0]}
+e_{\sssize [2,1]}\,e_{\sssize [2,0]}
+e_{\sssize [2,1]}\,L^2\,+
\\
+\,e_{\sssize [2,0]}\,e_{\sssize [0,1]}\,L^2
-e_{\sssize [1,1]}\,e_{\sssize [1,0]}\,L^2
+e_{\sssize [0,1]}\,L^4.
\endgathered
\mytag{5.20}
$$
Transforming \mythetag{5.20} with the use of the formulas \mythetag{5.9}, 
\mythetag{5.4}, \mythetag{4.19}, \mythetag{4.12}, we get
$$
\hskip -2em
\gathered
-4\,e_{\sssize [0,1]}\,e_{\sssize [1,1]}^2
+4\,e_{\sssize [1,1]}\,e_{\sssize [1,2]}
-2\,e_{\sssize [1,1]}\,e_{\sssize [1,0]}^3
+6\,e_{\sssize [2,1]}\,e_{\sssize [1,0]}^2\,+
\\
+\,2\,e_{\sssize [2,1]}\,e_{\sssize [0,1]}^2
+e_{\sssize [0,1]}^3\,e_{\sssize [1,0]}^2
-e_{\sssize [0,1]}^5+2\,e_{\sssize [2,1]}\,L^2\,-
\\
-\,2\,e_{\sssize [1,1]}\,e_{\sssize [1,0]}\,L^2
+4\,e_{\sssize [0,1]}\,e_{\sssize [1,0]}^2\,L^2
+7\,e_{\sssize [0,1]}^3\,L^2-4\,e_{\sssize [0,1]}\,L^4=0.
\endgathered
\mytag{5.21}
$$
Then we apply \mythetag{5.12} to \mythetag{5.21} in order to eliminate the
term with the square $e_{\sssize [1,1]}^2$:
$$
\pagebreak
\hskip -2em
\gathered
4\,e_{\sssize [1,1]}\,e_{\sssize [1,2]}
-2\,e_{\sssize [1,1]}\,e_{\sssize [1,0]}^3
+6\,e_{\sssize [2,1]}\,e_{\sssize [1,0]}^2\,+
\\
+\,2\,e_{\sssize [2,1]}\,e_{\sssize [0,1]}^2
-e_{\sssize [0,1]}^3\,e_{\sssize [1,0]}^2
+e_{\sssize [0,1]}\,e_{\sssize [1,0]}^4
+2\,e_{\sssize [2,1]}\,L^2\,-
\\
-\,2\,e_{\sssize [1,1]}\,e_{\sssize [1,0]}\,L^2
+2\,e_{\sssize [0,1]}\,e_{\sssize [1,0]}^2\,L^2
+e_{\sssize [0,1]}^3\,L^2
-3\,e_{\sssize [0,1]}\,L^4=0.
\endgathered
\mytag{5.22}
$$
The equation \mythetag{5.22} is the eighth factor equation derived from 
the cuboid equations \mythetag{1.1}. It is similar to the equation
\mythetag{5.17}.\par
\head
6. Concluding remarks.
\endhead
     One can continue deriving factor equation more and more. In order
to reasonably terminate this process we need some theoretical considerations.
Note that the left hand sides of the factor equations are polynomials from the 
ring
$$
\hskip -2em
\Bbb Q[e_{\sssize [1,0]},e_{\sssize [2,0]},e_{\sssize [3,0]},
e_{\sssize [0,1]},e_{\sssize [0,2]},e_{\sssize [0,3]},
e_{\sssize [2,1]},e_{\sssize [1,1]},e_{\sssize [1,2]},L],
\mytag{6.1}
$$
where $e_{\sssize [1,0]}$, $e_{\sssize [2,0]}$, $e_{\sssize [3,0]}$,
$e_{\sssize [0,1]}$, $e_{\sssize [0,2]}$, $e_{\sssize [0,3]}$, $e_{\sssize [2,1]}$,
$e_{\sssize [1,1]}$, $e_{\sssize [1,2]}$, and $L$ are treated as independent 
variables. If we continue deriving factor equation endlessly, their left hand 
sides would generate a certain ideal $J$ in the ring \mythetag{6.1}. By means 
of the formulas \mythetag{3.2}, \mythetag{3.3}, and \mythetag{3.4} this ideal 
$J$ is mapped onto some certain ideal $I_{\text{sym}}$ of the ring of multisymmetric 
polynomials $\Sym\!\Bbb Q[x_1,x_2,x_2,d_1,d_2,d_3,L]$. The ideal $I_{\text{sym}}$ 
is produced as the intersection 
$$
\hskip -2em
I_{\text{sym}}=I\cap\Sym\!\Bbb Q[x_1,x_2,x_2,d_1,d_2,d_3,L],
\mytag{6.2}
$$
where $I$ is the ideal of the polynomial ring $\Bbb Q[x_1,x_2,x_2,d_1,d_2,d_3,L]$
generated by the left hand sides of the cuboid equations \mythetag{1.1}.
Calculating the intersection \mythetag{6.2} is an algorithmically solvable 
computational problem. It will be considered in a separate paper.

\enddocument
\end